\documentclass[a4paper]{amsart}
\usepackage{amssymb}
\usepackage{amsmath}
\usepackage{amsthm}
\DeclareMathOperator*{\vol}{vol}
\begin{document}
\thanks{The first author was
 partially supported by  PSC CUNY Research Award,
 No. 60007-33-34, and George Shuster Fellowship at Lehman College}

\newtheorem{introtheorem}{Theorem}
\renewcommand{\theintrotheorem}{\Alph{introtheorem}}
\newtheorem{theorem}{Theorem }[section]
\newtheorem{lemma}[theorem]{Lemma}
\newtheorem{corollary}[theorem]{Corollary}
\newtheorem{proposition}[theorem]{Proposition}
\theoremstyle{definition}
\newtheorem{definition}[theorem]{Definition}
\newtheorem{example}[theorem]{Example}
\newtheorem{remark}[theorem]{Remark}

\renewcommand{\labelenumi}{(\roman{enumi})} 
\def\theenumi{\roman{enumi}}

\numberwithin{equation}{section}

\def \g {{\gamma}}
\def \G {{\Gamma}}
\def \l {{\lambda}}
\def \a {{\alpha}}
\def \b {{\beta}}
\def \f {{\phi}}
\def \r {{\rho}}
\def \R {{\mathbb R}}
\def \H {{\mathbb H}}
\def \N {{\mathbb N}}
\def \C {{\mathbb C}}
\def \Z {{\mathbb Z}}
\def \F {{\Phi}}
\def \Q {{\mathbb Q}}
\def \e {{\epsilon }}
\def \ev {{\vec\epsilon}}
\def \ov {{\vec{0}}}
\def \GinfmodG {{\Gamma_{\!\!\infty}\!\!\setminus\!\Gamma}}
\def \GmodH {{\Gamma\backslash\H}}
\def \sl  {{\hbox{SL}_2( {\mathbb R})} }
\def \psl  {{\hbox{PSL}_2( {\mathbb R})} }
\def \L  {{\hbox{L}^2}}

\newcommand{\norm}[1]{\left\lVert #1 \right\rVert}
\newcommand{\abs}[1]{\left\lvert #1 \right\rvert}
\newcommand{\modsym}[2]{\left \langle #1,#2 \right\rangle}
\newcommand{\inprod}[2]{\left \langle #1,#2 \right\rangle}
\newcommand{\Nz}[1]{\left\lVert #1 \right\rVert_z}
\newcommand{\tr}[1]{\operatorname{tr}\left( #1 \right)}

\title[Poincar\'e pairing, Riemann surfaces]{The distribution of values of the Poincar\'e pairing for hyperbolic Riemann surfaces}
\author{Yiannis N. Petridis}
\address{Department of Mathematics and Computer Science\\
City University of New York, Lehman College\\
250 Bedford Park Boulevard West,
Bronx\\NY 10468-1589}
\email{petridis@comet.lehman.cuny.edu}
\author{Morten S. Risager}
\address{Department of Mathematical Sciences\\ University of Aarhus\\ Ny Munkegade Building 530\\ 8000 {Aa}rhus, Denmark}
\email{risager@imf.au.dk}
\date{\today}
\begin{abstract}
 For a cocompact group of $\sl$ we fix a non-zero harmonic $1$-form
 $\alpha$. We normalize and order the values of the Poincar\'e pairing $\modsym{\gamma}{\a}$ according to the length of the corresponding closed geodesic $l(\gamma )$. We prove that these normalized values have a Gaussian distribution.

\end{abstract}
\maketitle
\section{Introduction}  
Let $X$ be a differentiable manifold.
We have the  pairing between homology and
cohomology: 
\begin{align*}H_1(X,\R)\times H^1_{\rm dR}(X,\mathbb R)\rightarrow \mathbb R\end{align*} and a projection
 $\phi:\Gamma =\pi_1(X)\rightarrow H_1(X,\Z)$. Let $\modsym{\cdot}{\cdot}$ be the
composition of the two maps:
$$\modsym{\gamma}{\a}=\int_{\phi (\gamma)}\a.$$
 We would like to study the distribution of the values of this for a fixed 
$1$-form $\a$. In previous work the authors \cite{petr, risager} have studied this problem for compact and finite volume hyperbolic surfaces. In both articles we found as limiting distribution the normal Gaussian distribution. However, the ordering of the group elements was not geometric: in \cite{petr} we ordered the group elements by realizing 
$\G=\pi_1(X)$ as a discrete subgroup of $\sl$, setting $\gamma=\left(\begin{array}{cc}a &b\\c&d\end{array}\right)$ and ordering $\gamma$ according to $c^2+d^2$. In \cite{risager} the matrix elements are ordered according to $(a^2+b^2)(c^2+d^2)$.
In both cases the ordering appears to be forced on us by the methods used: Eisenstein series associated with the problem in which the group elements are naturally summed in the above-mentioned fashion. 
There is a more natural geometric ordering. To every conjugacy class $\{\gamma\}$ corresponds
a unique closed oriented geodesic of length $l(\gamma)$. Let $\pi (x)=\#\{\{\gamma\}|
 \textrm{$\g$ prime},  l(\gamma)\le x\}$.
The prime number theorem for closed geodesics states that \begin{equation}\pi (x)\sim e^x/x\end{equation} 
as $ x\rightarrow \infty$ and can be proved using the Selberg trace formula (\cite{huber, buser}). In this article we consider the distribution of the values of the Poincar\'e pairing
where we order the elements of $\Gamma$ according to the lengths $l(\gamma)$.
\begin{theorem}\label{maintheorem}Let 
$$[\gamma, \a]=\sqrt{\frac{\vol(X)}{2\norm{\a}^2l(\gamma)}}\int_{\phi(\gamma )}\a.$$
 Then
 $$\frac{\#\left\{\gamma\in \pi_1(X)| [\gamma, \a]\in [a, b], l(\gamma)\le x\right\}}{\#\{\gamma\in\pi_1(X)|l(\gamma )\le x\}}\rightarrow\frac{1}{\sqrt{2\pi}}\int_a^be^{-t^2/2}\, dt $$
 as $x\rightarrow\infty$.
 \end{theorem}
 To prove this theorem  we use the method of moments, see
 \cite{loeve}. 

 Our approach is quite traditional: the Selberg trace formula, via the Selberg zeta function and its derivatives
 in character varieties. The geometric side gives us estimates for sums of $\modsym{\gamma}{\a}^n$ 
 ordered according to $l(\gamma)$. The spectral side involves the spectrum of the Laplace operator, as encoded in its
 resolvent $R(s)=(\Delta+s(1-s))^{-1}$. To extract information from
 the spectral side, we use perturbations of the resolvent.
 
 We follow the spirit of our previous work, which was motivated  by the question of finding the distribution and the moments of modular symbols. For this purpose Goldfeld \cite{goldfeld1, goldfeld2} introduced Eisenstein series associated with modular symbols. In \cite{risager} the second author introduced and studied the properties of hyperbolic Eisenstein series associated with modular symbols. Our current work uses the Selberg zeta function and its derivatives in various directions in character varieties. Such perturbations were first studied by Fay \cite{fay}.  

\begin{remark}
The study of the first moments
$$\pi (x, \a)= \sum_{\substack{\gamma\in \pi_1(X)\\l(\gamma)\le x}}\int_{\phi(\gamma )}\a$$
was initiated by Zelditch \cite{zelditch1, zelditch2} in relation to Bowen's equidistribution theorem for closed geodesics.
He proved bounds of the form $\pi (x, \a)=o(x/\ln (x))$ for 
$\a$ an automorphic form (perpendicular to the constants). He treated also finite-area hyperbolic surfaces. The technique in \cite{zelditch1, zelditch2} is the trace formula for the composition of two operators: standard convolution with a point-pair invariant, followed by multiplication by $\a$. Our work shows how (in principle) one could get asymptotics for $\pi (x, \a)$
depending on the Laurent series of the resolvent at $s=1$.

\end{remark}
\begin{remark}
Estimating the number of closed geodesics (periods of the geodesic
flow) with length less than $x$ and certain constraints on its periods can be done using Lalley's theorem \cite{lalley}, see also \cite{ledrappier}, \cite{sharp}. The constraints discussed in these articles restrict the periods to lie in a compact interval, say, $[a, b]$. In our theorem  we roughly restrict the periods
to lie in intervals $[a\sqrt{l(\gamma )}, b\sqrt{l(\gamma )}]$. In \cite{sharp2} a Gaussian law  similar to ours is stated, as a consequence of \cite{lalley} and \cite{ledrappier}. Sharp \cite{sharp2} also gets a local limit theorem. The method in all these papers is to use the thermodynamic formalism and as a consequence the results apply to variable negative curvature.
We stick to a more classical approach using the Selberg trace formula and get more explicit results.
\end{remark}
\section{The Selberg Zeta Function}
Let $\G$ be a discrete subgroup of $\psl$ with compact quotient $X=\GmodH$. Here $\H$ is the upper half-plane. The Selberg zeta function is defined as
\begin{equation}\label{selbergzeta}
Z(s, \chi )=\prod_{\{\gamma_0\}}\prod_{k=0}^{\infty}(1-\chi (\gamma_0)N(\gamma_0)^{-s-k}),\end{equation}
where $\chi$ is a unitary character of $\Gamma$ and the product is over primitive conjugacy classes, i.e. $\gamma_0$
is not a power of another element in $\Gamma$. The norm $N(\gamma)$ is
defined as follows. We conjugate (in $\sl$) the hyperbolic matrix
$\gamma$ to $\left(\begin{array}{cc}m&0\\0&m^{-1}\end{array}\right)$
with $m>1$. Then $N(\gamma)=m^2$. There is a simple relation between $N(\gamma)$ 
and $l(\gamma)$: $l(\gamma)=\log(N(\gamma))$. The product in (\ref{selbergzeta}) converges absolutely for $\Re (s)>1$.
We shall use the following convention: if in a product or sum the group elements carry a subscript $0$, then it extends over primitive elements/conjugacy classes. If not, then it extends over all group elements/conjugacy classes. 
For our purposes we need in fact a
family of characters 
$$\begin{array}{lccc}\chi(\cdot,\e):&\Gamma&\rightarrow& S^1\\
                   &\g&\mapsto &\exp(-i\e\modsym{\g}{\alpha}).
\end{array}
$$
We will denote the corresponding Selberg zeta functions by $Z(s, \e)$.
As in \cite{hejhal} ,Prop. 3.5, p.~56,
4.2, p.~67, we have the following splitting of $Z'(s, \e)/Z(s, \e)$:
\begin{equation}\label{zetasplitting}
\frac{Z'(s,\e)}{Z(s, \e)}=\sum_{\{\gamma_0\}}\frac{\chi (\gamma_0, \e)\ln (N(\gamma_0))}{N(\gamma_0 )^s}+A_1(s, \e)+A_2(s, \e),
\end{equation}
where 
\begin{equation}\label{Karl}
A_1(s, \e)=\sum_{\{\gamma\}}\frac{\chi (\gamma, \e)\ln (N(\gamma_0))}{N(\gamma )^s[N(\gamma)-1]}
\end{equation}
and
\begin{equation}\label{Sigrid}
A_2(s, \e)=\sum_{\{\gamma_0\}}\ln (N(\gamma_0))\frac{\chi (\gamma_0, \e)^2N(\gamma_0)^{-2s}}{1-\chi (\gamma_0, \e)N(\gamma_0)^{-s}}.
\end{equation}
We have the following lemma based on \cite{eichler}:
\begin{lemma}\label{eichlerbound}
\begin{equation*}
\modsym{\g}{\a} =\int_{\phi (\gamma )}\a=O( \ln N(\gamma)),
\end{equation*}
where the implied constant depends on the group only.
\end{lemma}
\begin{proof}
We let $F\subseteq \H$ be the normal polygon with $i$ in its interior. We may pick
generators $\g_1,\ldots\g_M$ such that $\g_j(F)$ is a neighboring
fundamental domain.  We may write every $\g\in\G$ as 
$$\g=\prod_{j=1}^{k(\g)} \g_{m_j},\qquad\textrm{ where $m_j\in\{1,\ldots,M\}$} $$
We let $C=\max\{\abs{\modsym{\g_j}{\a}}\vert j=1\ldots M\}$. Clearly
$\abs{\modsym{\g}{\a}}\leq Ck(\g)$. We now quote \cite{eichler}, Satz
1,  to conclude the existence of constants $a_\G, b_\G$ such
that \begin{equation}\label{eichlerbound2}k(\g)\leq a_\G\log(a^2+b^2+c^2+d^2)+b_\G.\end{equation}
Here $a,b,c,d$ are the entries of  $\g$.  

 Pick $B>0$ such that $B\geq d(i,x)$ for all $x\in
F$. Here $d(\cdot,\cdot)$ is the hyperbolic distance in $\H$. Such $B$ exists since
$X$ is compact.
Since $l(\g)=\inf_{x\in
  F}d(x,\g x)$ we conclude from 
$$d(i,\g i)\leq d(i,x)+d(x,\g x)+d(\g x,\g i)=d(x,\g x)+2d(i,x)\leq
d(x,\g x)+2B$$ that $d(i,\g i)\leq l(\g)+2B$. One may check that 
$$\cosh d(i,\g i)=\left(1+\frac{\abs{i-\g i}^2}{2\Im (i)\Im(\g
    i)}\right)=\frac{a^2+b^2+c^2+d^2}{2}.$$
Hence
\begin{align*}
\log(a^2+b^2+c^2+d^2)&=\log(2\cosh (d(i,\g i)))\\
&\leq \log(2
(\cosh( l(\g)+2B)))\\
&=\log (e^{\log N(\g)+2B}+e^{-\log N(\g)-2B})=O(\log(N(\g)))
\end{align*}
which finishes the proof.
\end{proof}
Using Lemma \ref{eichlerbound} and \cite{hejhal}, Prop. 2.5, we easily see that the series $A_1(s, \e)$ converges absolutely for $\Re (s)>0$, that
$A_2(s, \e)$ converges absolutely for $\Re (s)>1/2$ and the same is true for  the derivatives in 
$\e$ of {\em all} order for $A_1(s, \e)$ and $A_2(s, \e)$.
To apply the standard methods of analytic number theory we study the series
\begin{equation}\label{smaglos}
E(s, \e)=\sum_{\{\gamma_0\}}\frac{\chi(\g_0,\e)\ln (N(\gamma_0))}{N(\gamma_0)^s}
\end{equation}
and its derivatives at $\e=0$
\begin{equation}\label{marcopolo}
E^{(n)}(s,0)=\sum_{\{\gamma_0\}}\frac{{(-i)^n}\modsym{\g_0}{\a}^n\ln (N(\gamma_0))}{N(\gamma_0)^s}.
\end{equation}
We need to know the pole order at $s=1$  and the leading term in the Laurent expansion of the derivatives in $\e$ of $Z'(s,\e)/Z(s,\e)$ at $\e=0$, and we need to control the derivatives as functions of $s$ on vertical lines for $\Re (s)>1/2$.
The series (\ref{marcopolo}) is similar to the Eisenstein series twisted by modular symbols
introduced by Goldfeld \cite{goldfeld1} and their holomorphic analogues introduced by Eichler \cite{eichler}.

\section{The automorphic Laplacian}
The beautiful connection between the Selberg zeta function and the spectrum of the automorphic Laplacian goes through  the Selberg trace formula. We shall only briefly touch upon this connection and refer to \cite{hejhal, buser, selberg2, venkov2} for further details.

We let $$\Delta=y^2\left(\frac{\partial^2}{\partial
    x^2}+\frac{\partial^2}{\partial y^2}\right)$$ be the Laplace
operator for the upper half-plane. We consider the space
$$\L(\GmodH,\overline \chi(\cdot,\epsilon))$$ of
$(\G,\overline \chi(\cdot,\epsilon))$-automorphic functions, i.e. functions
$f:\H\to\C$ where $$f(\g z)=\overline{\chi}(\g,\e)f(z),$$ and
$$\int_{\GmodH}\abs{f(z)}^2d\mu(z).$$  Here $d\mu(z)=y^{-2}dxdy$ is
the invariant Riemannian measure on $\H$ derived from the Poincar\'e
metric $ds^ 2=y^{-2}(dx^2+dy^2)$. We shall denote by $\norm{\cdot}$
the usual norm in the Hilbert space $\L(\GmodH,\overline
\chi(\cdot,\epsilon))$. The automorphic Laplacian $\tilde L(\e)$ is
the closure of the operator acting on smooth functions in
$\L(\GmodH,\overline \chi(\cdot,\epsilon))$ by $\Delta f$. The  spectrum of $\tilde{L}(\e)$ is discrete and $-\tilde L(\e)$ is nonnegative with eigenvalues
$$0\leq \l_0(\e)\leq \l_1(\e)\leq\ldots ,$$ 
satisfying $\lim_{n\to\infty}\l_n(\e)=\infty$,  \begin{equation}\sum_{n=1}^\infty \l_n(\e)^{-2}<\infty,\end{equation}  and 
\begin{equation}
\#\left\{\l_n(\e)\leq \l \right\}\sim\frac{\vol{(\GmodH)}}{4\pi}\l\textrm{ as }\l\to\infty.
\end{equation} The corresponding eigenfunctions  $\{\phi_n(\e)\}_{n=0}^ \infty$ may be chosen such that they form a complete orthonormal family.  The first eigenvalue is zero if and only if $\e=0$ and in this case it is a simple eigenvalue. 
We consider two additional spectral problems on $\GmodH$: the Dirichlet problem and the Neumann (free boundary)
problem. We denote their corresponding eigenvalues by $l_n$ and $\mu_n$ respectively. Using the Rayleigh quotient we easily get: $\mu_n \le \l_n(\e)\le l_n$ for  all $\e$. This implies for the spectral counting functions  
\begin{equation}\label{counter}
\frac{\vol{(\GmodH)}}{4\pi}\l\sim N_D(\l)\le N_{\e}(\l)\le N_N(\l)\sim\frac{\vol{(\GmodH)}}{4\pi}\l\textrm{ as }\l\to\infty.
\end{equation}
The resolvent $\tilde R(s,\e)=(\tilde L(\e)+s(1-s))^{-1}$, defined off the spectrum of $\tilde L(\e)$, is a Hilbert-Schmidt $\L(\GmodH,\overline \chi(\cdot,\epsilon))$ operator. It is holomorphic in $s$, and its operator norm may be bounded as 
\begin{equation}\label{calculustestsareboring}
\norm{\tilde R(s,\e)}_\infty\leq \frac{1}{\hbox{dist}({s(s-1),\hbox{spec}(\tilde L(\e))})}\leq \frac{1}{\abs{t}(2\sigma-1)},
\end{equation} where $s=\sigma+it$, $\sigma>1/2$.
We recall that for a compact operator, $T$, the singular values, $\{\beta_k\}_{k=0}^\infty$, are defined as the square roots of the eigenvalues of $T^*T$, i.e. 
$$T^*T f=\sum_{k=1}^\infty \beta_k^2(f,\psi_k)\psi_k,\qquad \beta_1\geq\beta_2\geq\ldots>0,$$ where $\{\psi_k\}_{k=0}^\infty$ forms a complete orthonormal family. When $T$ is symmetric the singular values are the absolute values of the eigenvalues of $T$. The $p$-norm is defined by
$$\norm{T}_p^p=\sum_{k=0}^\infty\beta_k^p.$$
When $p=2$ this is the Hilbert-Schmidt norm and when $p=1$ this is the trace norm.

We shall need a bound on the Hilbert-Schmidt norm of the resolvent. We let $s=\sigma+it$.
\begin{lemma}For fixed $\sigma> 1/2$ the Hilbert-Schmidt norm of the resolvent is bounded as $\abs{t}\to\infty$. More precisely \begin{equation}\label{resolventbound} 
\norm{\tilde R(s,\e)}_2\leq O((2\sigma-1)^{-1})\end{equation} as $t\to\infty$. The involved constant depends only on the group $\G$.\end{lemma}
\begin{proof}
Let $\l_n(\e)=1/4+r_n(\e)^2$ with $r_n(\e)\in\R_+\cup i\R_+$. Then $s(1-s)-\l_n(\e)=-((s-1/2)^2+r_n(\e)^2)$. Hence 
\begin{align}
\norm{\tilde R(s,\e)}_2^2=&\sum_{n=0}^\infty\frac{1}{\abs{(s-1/2)^2+r_n(\e)^2}^2}\\\allowdisplaybreaks
\nonumber =&\sum_{n=0}^\infty\frac{1}{\abs{(\sigma-1/2)+i(t+r_n(\e))}^2\abs{(\sigma-1/2)+i(t-r_n(\e))}^2} 
\end{align}
For $r_n(\e)\in iR_+$ all individual summands is less than $t^{-4}$ and, since there can only be finitely many small eigenvalues , the sum over small eigenvalues is $O_\e(t^{-4})$.

For $0\leq r_n(\e) \leq 2t$ the individual summands is less than $((\sigma-1/2)t)^{-2}$ and (\ref{counter}) says that there is $O(t^2)$ elements in the sum. Hence
\begin{equation*}
\sum_{0\leq r_n(\e)\leq 2t}\frac{1}{\abs{(s-1/2)^2+r_n(\e)^2}^2}=O((\sigma-1/2)^{-2}).
\end{equation*} 

When $r_n(\e)>2t$ we have $r_n(\e)-t>r_n(\e)/2$ and hence the individual terms in the sum may be bounded by $r_n(\e)^{-4}$. We therefore have 
\begin{equation*}
  \sum_{r_n(\e)>2t}\frac{1}{\abs{(s-1/2)^2+r_n(\e)^2}^2}\leq 4 \sum_{r_n(\e)>2t}\frac{1}{r_n(\e)^4}=O(t^{-2})
\end{equation*} 
which finishes the proof. (We have again used  (\ref{counter}) for the last equality.)
\end{proof}
We fix $z_0\in\H$ and introduce unitary operators
\begin{equation}\begin{array}{rccc}
U(\e):&\L(\GmodH)&\to& \L(\GmodH,\overline \chi(\cdot,\epsilon))\\
&f&\mapsto&\exp\left(i\e\int_{z_0}^z\alpha\right)f(z). 
\end{array}
\end{equation}
We then define
\begin{eqnarray}
L(\e)&=&U^{-1}(\e)\tilde L(\e)U(\e)\\
R(s,\e)&=&U^{-1}(\e)\tilde R(s,\e)U(\e). 
\end{eqnarray}
This ensures that $L(\e)$ and $R(s,\e)$ act on the fixed space $\L(\GmodH)$. It is then easy to verify that 
\begin{align}\label{muffinstogo}L(\e)h=\Delta h +2i\e\modsym{dh}{\alpha}&-i\e\delta(\alpha)h-\e^ 2\modsym{\a}{\a}\\
\label{coffeetogo}(L(\e)+s(1-s))R(s,\e)=&R(s,\e)(L(\e)+s(1-s))=I.
\end{align}
Here \begin{eqnarray*}\modsym{f_1dz+f_2d\overline z}{g_1dz+g_2d\overline z}&=&2y^ 2(f_1\overline{g_1}+f_2\overline{g_2})\\
\delta(pdx+qdy)&=&-y^2(p_x+q_y).
\end{eqnarray*}
We notice that $\delta(\alpha)=0$. We notice also that 
\begin{align}
L^{(1)}(\e)h&=2i\modsym{dh}{\a}-2\e\modsym{\a}{\a},\\ \allowdisplaybreaks
L^{(2)}(\e)h&=-2\modsym{\a}{\a},\\\allowdisplaybreaks
L^{(i)}(\e)h&=0,\quad\textrm{ when }i\geq 3.
\end{align}
(We use superscript $(n)$ to denote the
$n$'th derivative in $\e$.)
Fix $\kappa>1$. The resolvent and the Selberg zeta function are
connected through the following identity, see \cite{hejhal}, Theorem 4.10, p.~72 ,:
\begin{align}
\nonumber\frac{1}{1-2s}&\frac{Z'}{Z}(s,\e)-\frac{1}{1-2\kappa}\frac{Z'}{Z}(\kappa,\e)=\\
\label{bigidentity}\sum_{n=0}^\infty&\left(\frac{1}{s(1-s)-\l_n(\e)}-\frac{1}{\kappa(1-\kappa)-\l_n(\e)}\right)-\frac{\vol(\GmodH)}{2\pi}\sum_{k=0}^\infty\left(\frac{1}{s+k}-\frac{1}{\kappa+k}\right).
\end{align}
 We note that by the Hilbert identity
\begin{equation}R(s,\e)-R(\kappa,\e)=-(s(1-s)-\kappa(1-\kappa))R(s,\e)R(\kappa,\e)
\end{equation}
the difference of resolvents is the product of two Hilbert-Schmidt operators. Hence, by the inequality \begin{equation}\label{ulighed} \norm{ST}_1\leq\norm{S}_2\norm{T}_2\end{equation} we conclude that $R(s,\e)-R(\kappa,\e)$ is of the trace class. We can therefore define the trace
\begin{equation}
\tr{R(s,\e)-R(\kappa,\e)}=\sum_{n=0}^\infty ((R(s,\e)-R(\kappa,\e))\psi_n,\psi_n).
\end{equation}
We note that this enables us to write (\ref{bigidentity}) in the form
\begin{equation}\label{niceidentity}
\frac{1}{1-2s}\frac{Z'}{Z}(s,\e)-\frac{1}{1-2\kappa}\frac{Z'}{Z}(\kappa,\e)=\tr{R(s,\e)-R(\kappa,\e)}+Q(s),
\end{equation} where $Q(s)$ is the last sum in (\ref{bigidentity}) (Note that $R(s,\e)-R(\kappa,\e)$ and $\tilde R(s,\e)-\tilde R(\kappa,\e)$ have the same trace). It is this identity that we shall study to obtain information about the analytic properties of (\ref{marcopolo}). 

We want to see how this behaves as $\abs{t}\to\infty$.
It is easy to see (\cite{hejhal}, p. 80 ) that if $s=\sigma+it$ then for fixed $\sigma>0$ we have $Q(s)=O(\log t)$.   From (\ref{resolventbound}) and (\ref{ulighed}) we conclude that when $\sigma > 1/2$ 
\begin{equation}\label{goon} \tr{R(s,\e)-R(\kappa,\e)}=O(t^2).\end{equation}
(By more careful estimates one may prove $O(t)$. See \cite{hejhal}, Eq. (4.8)) Note that when $\e=0$ this proves that $E(z,s)$ grows at most like $O(t^3)$ on
vertical lines $\sigma>1/2$ (See (\ref{zetasplitting})). We wish to prove something
similar for $E^{(n)}(s,\e)$. To obtain this we  differentiate (\ref{coffeetogo}) in $\e$ and get 
\begin{align}\label{blah}R^{(n)}(s,\e)=&-\sum_{i=0}^{n-1}\binom{n}{i}R^{(i)}(s,\e)L^{(n-i)}(\e)R(s,\e)\\
\nonumber=&-\sum_{i=1}^{n}\binom{n}{i}R(s,\e)L^{(i)}(\e)R^{(n-i)}(s,\e).
\end{align}
This is naturally the second Neumann series for the resolvent.
We use $L^{(i)}(\e)=0$ for $i\geq 3$, which  reduces (\ref{blah}) to 
\begin{align}\label{comehomelaptop}R^{(1)}(s,\e)&=-R(s,\e)L^{(1)}(\e)R(s,\e)
 \\\nonumber R^{(n)}(s,\e)&=-\left(\binom{n}{1}R(s,\e)L^{(1)}(\e)R^{(n-1)}(s,\e)+\binom{n}{2}R(s,\e)L^{(2)}(\e)R^{(n-2)}(s,\e) \right)
\end{align}
We need the following lemma:
\begin{lemma}\label{norms}The operators $R(s,\e)L^{(i)}(\e)$,  $L^{(i)}(\e)R(s,\e)$ are bounded and their norms grow at most polynomially for $s$ on a fixed line $\Re(s)=\sigma>1/2$. More precisely we have
\begin{align}
\label{firstnorms}\norm{R(s,\e)L^{(1)}(\e)}_\infty&=O(\abs{t}),\qquad &\norm{L^{(1)}(\e)R(s,\e)}_\infty&=O(\abs{t}),\\
\label{secondnorms}\norm{R(s,\e)L^{(2)}(\e)}_\infty&=O\left(\abs{t}^{-1}\right),\qquad &\norm{L^{(2)}(\e)R(s,\e)}_\infty&=O\left(\abs{t}^{-1}\right).
\end{align} 
\end{lemma}
\begin{proof}
The claim in (\ref{secondnorms}) follows easily from (\ref{calculustestsareboring}) and the fact that $L^{(2)}(\e)$ is bounded since it is just a multiplication operator on a compact set. The first claim in (\ref{firstnorms}) follows from the second since 
$$\norm{R(s,\e)L^{(1)}(\e)}_\infty=\norm{\left(R(s,\e)L^{(1)}(\e)\right)^{*}}_\infty=\norm{L^{(1)}(\e)R(\overline{s},\e)}_\infty.$$ 
To prove the remaining case we  use Sobolev $s$-norms, $\norm{\cdot}_{H^s}$ and the fact that for any second order elliptic operator $P$ there exist a $c'$ such that
$$\norm{u}_{H^2}\leq c'(\norm{u}+\norm{Pu}).$$ 
We shall use $P=L(\e)$. Hence \begin{align*}
\norm{L^{(1)}(R(s,\e))u}\leq &c\norm{R(s,\e)u}_{H^1}\\
                         \leq  &c\norm{R(s,\e)u}_{H^2}\\
                         \leq    &c'( \norm{R(s,\e)u}+\norm{L(\e)R(s,\e)u})\\
= &c'( \norm{R(s,\e)u}+\norm{s(1-s)R(s,\e)u+u}) 
\end{align*}
where we have used (\ref{coffeetogo}). The result now follows from (\ref{calculustestsareboring}).
\end{proof}
\begin{theorem}\label{growth}For fixed $\sigma=\Re(s)$ the difference $$R^{(n)}(s,\e)-R^{(n)}(\kappa,\e)$$ is of the trace class and the trace grows at most polynomially in $\abs{t}=\abs{\Im(s)}$. For $1/2<\sigma \leq 1$ we have
$$\tr{R^{(n)}(s,\e)-R^{(n)}(\kappa,\e)}=O(\abs{t}^{n+2+\varepsilon}).$$
\end{theorem}
\begin{proof} Clearly $\abs{\tr{\cdot}}\leq \norm{\cdot}_1$. From (\ref{comehomelaptop}) it is clear that $R^{(n)}(s,\e)$ is a linear combination of terms of the form $R(s,\e)\prod_{i=1}^m(L^{k_i}(\e)R(s,\e))$ where $k_i\in\{1,2\}.$ It is also clear that one of the terms is a constant times $R(s,\e)(L^{(1)}(\e)R(s,\e))^n$ and that none of the products have more than $n$ terms. We shall prove by induction that \begin{equation*}\norm{R(s,\e)\prod_{i=1}^m(L^{k_i}(\e)R(s,\e))-R(\kappa,\e)\prod_{i=1}^m(L^{k_i}(\e)R(\kappa,\e))}_1=O(\abs{t}^{2+n})\end{equation*}
for $m\le n$.  
The case $n=0$ is (\ref{goon}). In the general case we add and subtract $R(s,\e)\prod_{i=1}^m(L^{k_i}(\e)R(\kappa,\e))$ and find
\begin{align*}
&\norm{R(s,\e)\prod_{i=1}^m(L^{k_i}(\e)R(s,\e))-R(\kappa,\e)\prod_{i=1}^m(L^{k_i}(\e)R(\kappa,\e))}_1\\
&\quad\leq  \norm{R(s,\e)\left(\prod_{i=1}^m(L^{k_i}(\e)R(s,\e))-\prod_{i=1}^m(L^{k_i}(\e)R(\kappa,\e))\right)}_1\\&\qquad+\norm{\left(R(s,\e)-R(\kappa,\e)\right)\prod_{i=1}^m(L^{k_i}(\e)R(\kappa,\e)) }_1\\
&\quad\leq  \norm{R(s,\e)L^{k_m}}_\infty\norm{R(s,\e)\prod_{i=1}^{m-1}(L^{k_i}(\e)R(s,\e))-R(\kappa,\e)\prod_{i=1}^{m-1}(L^{k_i}(\e)R(\kappa,\e))}_1\\&\qquad+\norm{R(s,\e)-R(\kappa,\e)}_1\norm{\prod_{i=1}^{m}(L^{k_i}(\e)R(\kappa,\e)) }_\infty.
\end{align*}
We quote Lemma \ref{norms} and use the induction hypothesis. This  completes the induction. We conclude that  
\begin{equation}\norm{R^{(n)}(s,\e)- R^{(n)}(\kappa,\e)}_1=O(\abs{t}^{2+n}).\end{equation} 
\end{proof}
Equation (\ref{zetasplitting}) together with the fact that all derivatives in $\e$ of  $A_1(s,\e)$ and  $A_2(s,\e)$ are absolutely convergent for $\Re(s)>1/2$ enables us to conclude from Theorem \ref{growth} that the function $E^{(n)}(s,\e)$ grows at most polynomially in $t$ on every fixed vertical line $\sigma>1/2$.
\section{The pole at $s=1$.}
In this section we identify the pole order and the leading term of
$E^{(n)}(s)=E^{(n)}(s,0)$ at $s=1$. (In $E^{(n)}(s,\e)$, $R(s,\e)$ and
$L^{(n)}(\e))$ we shall often omit $0$ from the notation when we
put $\e=0$).  We note that $E(s)$ has a first order pole with residue 1 at this point as is easily seen from (\ref{bigidentity}). Since
\begin{equation}\label{meneldor}\left.\frac{d^n}{d\e^n}\left(\frac{1}{1\! -\! 2s}\frac{Z'}{Z}(s,\e)-\frac{1}{1\! -\! 2\kappa}\frac{Z'}{Z}(\kappa,\e)\right)\right\vert_{\e=0}\!\!\!\!=\!\tr{R^{(n)}(s)-R^{(n)}(\kappa)}
\end{equation} equation (\ref{comehomelaptop}) and the fact that $R^{(n)}(s)$ is holomorphic in $\Re(s)>1$ enables us to conclude that the left-hand side is holomorphic in $\Re(s)>1$.

We recall that close to $s=1$ \begin{equation}\label{resolventexpansion}R(s)=\sum_{i=-1}^\infty R_i (s-1)^i,\qquad R_{-1}=-P_0\end{equation}  and that $R(s)-R_{-1}(s-1)^{-1}$ is holomorphic in $\Re(s)>h$. Here  $h=\Re(s_1)$, $s_1(1-s_1)=\l_1$ is the first small eigenvalue, and $P_0f=\langle f,\vol(\GmodH)^{-1/2}\rangle\vol{(\GmodH)}^{-1/2}$ is the projection of $f$ to the zero eigenspace.

To understand the meromorphic structure of (\ref{meneldor}) at $s=1$, we must understand the meromorphic structure of $R^{(n)}(s)$. The crucial observation is that 
\begin{equation}\label{everythingdies}
L^{(1)} P_0=0, \qquad P_0 L^{(1)} =0.\end{equation}
The first equality follows from the fact that $L^{(1)}$ is a differentiation operator while $P_0$ projects to the constants. The second equality follows from the first by using that both operators are selfadjoint.
Using this we can now prove:
\begin{lemma}For $n\geq 0$, $R^{(n)}(s)$  has a singularity at $s=1$ of at most order $[n/2]+1$. If $n=2m$ the singularity is of order $m+1$ and the $(m+1)$-term in the expansion around $s=1$ is $$\frac{(2m)!}{2^ m}(-1)^m R_{-1}(L^{(2)} R_{-1})^m.$$
\end{lemma}
\begin{proof}
The claim about $n=0$ is contained in (\ref{resolventexpansion}). The proof is induction in $n$ but for clarity we do $n=1,2$ by hand. 
Using (\ref{comehomelaptop}) and (\ref{resolventexpansion}) we see that $R^{(1)}(s
,0)$ has a pole of  at most  second order and that the singular part is
$$R_{-1}L^{(1)} R_{-1}(s-1)^{-2}+(R_{0}L^{(1)} R_{-1}+R_{-1}L^{(1)} R_{0})(s-1)^{-1}.$$
But this is zero by (\ref{everythingdies}). Hence $R^{(1)}(s)$ is regular at $s=1$.

For $n=2$ we note that by (\ref{resolventexpansion}) and (\ref{everythingdies}) $R(s)L^{(1)} $ is regular. Hence we find that the first term in 
$$R^{(2)}(s)=-\left(\binom{2}{1}R(s)L^{(1)} R^{(1)}(s)+\binom{2}{2}R(s)L^{(2)} R(s)\right)$$ is regular. By (\ref{resolventexpansion}) the second term is of at most order two with leading term $-\binom{2}{2}R_{-1}L^{(2)} R_{-1}$.

In the general case we note that, since $R(s)L^{(1)} $ is regular the first term in 
$$R^{(n)}(s)=-\left(\binom{n}{1}R(s)L^{(1)} R^{(n-1)}(s)+\binom{n}{2}R(s)L^{(2)}R^{(n-2)}(s) \right)$$
has at most a pole of the same order as $R^{(n-1)}(s)$. By (\ref{resolventexpansion}) the second term has a pole of order at  most one more than that of $R^{(n-2)}(s)$. The claim about the order of the pole follows. The claim about the leading term also follows once we note that  
 $$\binom{2m}{2}\ldots\binom{6}{2}\binom{4}{2}\binom{2}{2}=\frac{(2m)!}{2^m}.$$
\end{proof}
Using this we can now prove the following theorem. Let $\norm{\alpha}^2=\int_{\GmodH}\modsym{\a}{\a} d\mu(z)$.
 \begin{theorem}\label{poles} The function $E^{(n)}(s)$  is identically zero for $n$ odd.  If $n=2m$ the function $E^{(n)}(z, s)$ has a pole  of order $m+1$ and the leading term is 
$$(-1)^m\frac{(2m)!\norm{\a}^{2m}}{\vol(\GmodH)^m}.$$
\end{theorem}
\begin{proof}
For $n$ odd the group elements $\gamma$ and $\gamma^{-1}$ contribute opposite values in (\ref{marcopolo}). When $n$ is even, we 
notice that by (\ref{niceidentity}) and (\ref{zetasplitting}) the leading term of  $E^{(n)}(s)$ is minus the leading term of $\tr{R^{(n)}(s)-(R^{(n)}(\kappa)}$. 
Using the above lemma, (\ref{zetasplitting}) and the fact that the derivatives of $A_1(s)$ and $A_2(s)$ are holomorphic in $\Re(s)>1/2$  we immediately get  the claim about the pole orders. The calculation of the leading term follows from the observation that 
\begin{align*}\sum_{k=0}^\infty (R_{-1}(L^{(2)} R_{-1})^m&\psi_k ,\psi_k )       =      - ((L^{(2)} R_{-1})^m\psi_0 ,\psi_0 )\\       =      &\frac{(-1)^{m+1}}{\vol{(\GmodH)}^m}\left(\int_{\GmodH}-2\modsym{\a}{\a}d\mu(z)\right)^m\end{align*}
\end{proof}

\section{Calculating the moments}
We are now ready to prove Theorem \ref{maintheorem}. The proof uses the method of asymptotical moments precisely as in \cite{petr,risager}. From Theorem \ref{poles}, Theorem \ref{growth} and Lemma \ref{eichlerbound} we may conclude, using a more or less standard contour integration argument (see  \cite{petr,risager} for details), that as $T\to\infty$
\begin{equation}\label{summatory}
\sum_{\substack{\{\g_0\}\\ N(\g_0)\leq T}}\modsym{\g_0}{\a}^n\log N(\g_0)=\left\{\begin{array}{ll}\displaystyle{\frac{(2m)!\norm{\a}^{2m}}{m!\vol{(\GmodH)}^{m}}}T(\log T)^m+O(T(\log T)^{m-1}),&n=2m,\\
0,&n=2m+1.
\end{array}\right.
\end{equation} 
Let now 
\begin{equation}
[\g,\a]=\sqrt{\frac{\vol({\GmodH})}{2\log(N(\g))\norm{\a}^2}}\modsym{\g}{\a}
\end{equation}
We then define the random variable $X_T$ with probability measure 
\begin{equation}
P(X_T\in [a,b])=\frac{
\#\{ \{ \g_0 \} | N(\g_0)\leq T,\quad[\g_0,\a] \in [a,b]\} }
{\#\{ \{ \g_0 \} | N(\g_0)\leq T\}}
\end{equation}
We want to calculate the asymptotical moments of these, i.e. find
\begin{equation}M_n(X_T)=\frac{1}{\#\{ \{ \g_0 \} | N(\g_0)\leq T\}}\sum_{\substack{\{ \g_0 \}\\ N(\g_0)\leq T}}[\g_0,\a]^n\end{equation} 
as $T\to\infty$. We note that by the prime number theorem for closed geodesics (or
(\ref{summatory})) the denominator is asymptotically $T/\log T$. By
partial summation we have 
\begin{align*}\sum_{N(\g_0)\leq T}&[\g_0,\a]^n=\frac{\vol(\GmodH)^{n/2}}{\norm{\a}^n2^{n/2}}\sum_{N(\g_0)\leq T}\modsym{\g_0}{\a}^n\log N(\g_0)\frac{1}{\log(N(\g_0))^{n/2+1}}\\
=&\frac{\vol(\GmodH)^{n/2}}{\norm{\a}^n2^{n/2}\log(T)^{n/2+1}}\sum_{N(\g_0)\leq T}\modsym{\g_0}{\a}^n\log N(\g_0)+O(\log(\log(T))),
\end{align*}
This may be evaluated by (\ref{summatory}). We find
$$M_n{(T)}\to \begin{cases}\displaystyle{\frac{(2m)!}{m!2^m}},&\textrm{if
      $n=2m$,}\\
0,&\textrm{otherwise.}
\end{cases}$$
We notice that the right-hand side coincides with the moments of the
Gaussian distribution. Hence by a classical result due to Fr\'echet
and Shohat (See \cite{loeve}, 11.4.C ) we may conclude that 
$$P(X_T\in
[a,b])\to\frac{1}{\sqrt{2\pi}}\int_a^b\exp\left(-\frac{x^2}{2}\right)dx\textrm{ as $T\to\infty$}.$$
This is almost Theorem \ref{maintheorem}. The only difference is that
here we are only counting the prime geodesics.  Let  $\mu>0$ be a lower bound for the lengths of the closed geodesics
on $\GmodH$. Define
\begin{align*}
\phi_n(T)=& \sum_{\substack{\{\gamma\}\\ N(\g)\leq T}}\modsym{\g}{\alpha}^n\\
\Pi_n(T)= & \sum_{\substack{\{\gamma_0\}\\ N(\g_0)\leq T}}\modsym{\g_0}{\alpha}^n\\
\end{align*}
i.e. the first sums over all conjugacy classes, while the second only
over the primitives. Since every conjugacy class may be written uniquely as the power
of a primitive one,  we use $\modsym{\g_0^m}{\a}=m \modsym{\g_0}{\a}$
and $N(\g_0^m)=N(\g_0)^m$ to
see that  
\begin{equation}\label{letsgetreal}\phi_n(T)=\Pi_n(T)+\sum_{m=2}^{[\log T/\log \mu]}m^n\Pi_n(\root m \of{T}).\end{equation}
The sum is clearly $O((\log T)^{n+1}\Pi_n(\root 2 \of{T}))$. Using
partial integration as above it is not difficult to see that 
\begin{equation}
\Pi_n(T)\sim \frac{1}{\log T}\sum_{\substack{\{\g_0\}\\ N(\g_0)\leq T}}\modsym{\g_0}{\a}^n\log N(\g_0).
\end{equation} 
Hence by (\ref{letsgetreal}) and the bound on the sum we have also
$\Pi_n(T)\sim\phi_n(T)$. Playing the same trick backwards we find that
(\ref{summatory}) is true also if we sum over  all conjugacy
classes and not only primitive ones. Doing the same argument as that
following (\ref{summatory}) we arrive at Theorem \ref{maintheorem}.

\begin{remark} We conclude by remarking that we could have proved a
  distribution result for the Poincar\'e  pairing between the homology
  classes and a fixed \emph{complex holomorphics} 1-form $f(z)dz$ instead of a
  harmonic 1-form $\alpha$. In this case the combinatorics involved
  would become more difficult as one needs to introduce an
  $n$-parameter family of character instead of $\chi(\cdot,\e)$. One
  finds that if we define
  $$[\g,f(z)dz]=\sqrt{\frac{\vol({\GmodH})}{2\log(N(\g))\norm{f(z)dz}^2}}\modsym{\g}{f(z)dz}$$
  and define the random variable $Y_T$ with probability measure 
\begin{equation}
P(Y_T\in R)=\frac{
\#\{ \{ \g \} | N(\g)\leq T,\quad[\g,f(z)dz] \in R\} }
{\#\{ \{ \g \} | N(\g)\leq T\}}
\end{equation}
where $R\subseteq \C$ is a rectangle, then 
$$P(Y_T\in R)\to\frac{1}{2\pi}\int_R\exp\left(-\frac{x^2+y^2}{2}\right)dxdy\textrm{ as $T\to\infty$}.$$
\end{remark}
{\bf Acknowledgments:}\newline{We would like to thank Erik Balslev, Fran\c cois Ledrappier, Peter Sarnak,  Alexei B. Venkov
and Steve Zelditch for valuable comments and 
suggestions. }

\end{document}